%
\magnification\magstep1
\baselineskip15pt
\newread\AUX\immediate\openin\AUX=\jobname.aux
\newcount\relFnno
\def\ref#1{\expandafter\edef\csname#1\endcsname}
\ifeof\AUX\immediate\write16{\jobname.aux gibt es nicht!}\else
\input \jobname.aux
\fi\immediate\closein\AUX



\def\ignore{\bgroup
\catcode`\;=0\catcode`\^^I=14\catcode`\^^J=14\catcode`\^^M=14
\catcode`\ =14\catcode`\!=14\catcode`\"=14\catcode`\#=14\catcode`\$=14
\catcode`\&=14\catcode`\'=14\catcode`\(=14\catcode`\)=14\catcode`\*=14
\catcode`+=14\catcode`\,=14\catcode`\-=14\catcode`\.=14\catcode`\/=14
\catcode`\0=14\catcode`\1=14\catcode`\2=14\catcode`\3=14\catcode`\4=14
\catcode`\5=14\catcode`\6=14\catcode`\7=14\catcode`\8=14\catcode`\9=14
\catcode`\:=14\catcode`\<=14\catcode`\==14\catcode`\>=14\catcode`\?=14
\catcode`\@=14\catcode`\A=14\catcode`\B=14\catcode`\C=14\catcode`\D=14
\catcode`\E=14\catcode`\F=14\catcode`\G=14\catcode`\H=14\catcode`\I=14
\catcode`\J=14\catcode`\K=14\catcode`\L=14\catcode`\M=14\catcode`\N=14
\catcode`\O=14\catcode`\P=14\catcode`\Q=14\catcode`\R=14\catcode`\S=14
\catcode`\T=14\catcode`\U=14\catcode`\V=14\catcode`\W=14\catcode`\X=14
\catcode`\Y=14\catcode`\Z=14\catcode`\[=14\catcode`\\=14\catcode`\]=14
\catcode`\^=14\catcode`\_=14\catcode`\`=14\catcode`\a=14\catcode`\b=14
\catcode`\c=14\catcode`\d=14\catcode`\e=14\catcode`\f=14\catcode`\g=14
\catcode`\h=14\catcode`\i=14\catcode`\j=14\catcode`\k=14\catcode`\l=14
\catcode`\m=14\catcode`\n=14\catcode`\o=14\catcode`\p=14\catcode`\q=14
\catcode`\r=14\catcode`\s=14\catcode`\t=14\catcode`\u=14\catcode`\v=14
\catcode`\w=14\catcode`\x=14\catcode`\y=14\catcode`\z=14\catcode`\{=14
\catcode`\|=14\catcode`\}=14\catcode`\~=14\catcode`\^^?=14
\Ignoriere}
\def\Ignoriere#1\;{\egroup}

\def\today{\number\day.~\ifcase\month\or
  Januar\or Februar\or M{\"a}rz\or April\or Mai\or Juni\or
  Juli\or August\or September\or Oktober\or November\or Dezember\fi
  \space\number\year}
\font\sevenex=cmex7
\scriptfont3=\sevenex
\font\fiveex=cmex10 scaled 500
\scriptscriptfont3=\fiveex
\def\A{{\bf A}}
\def\G{{\bf G}}
\def\P{{\bf P}}

\def\CS{{\widetilde C}}

\def\phi{\varphi}
\def\epsilon{\varepsilon}
\def\theta{\vartheta}
\def\uauf{\lower1.7pt\hbox to 3pt{%
\vbox{\offinterlineskip
\hbox{\vbox to 8.5pt{\leaders\vrule width0.2pt\vfill}%
\kern-.3pt\hbox{\lams\char"76}\kern-0.3pt%
$\raise1pt\hbox{\lams\char"76}$}}\hfil}}

\def\title#1{\par
{\baselineskip1.5\baselineskip\rightskip0pt plus 5truecm
\leavevmode\vskip0truecm\noindent\font\BF=cmbx10 scaled \magstep2\BF #1\par}
\vskip1truecm
\leftline{\font\CSC=cmcsc10\CSC Friedrich Knop}
\leftline{Department of Mathematics, Rutgers University, New Brunswick NJ
08903, USA}
\leftline{knop@math.rutgers.edu}
\vskip1truecm
\par}

\def\cite#1{\expandafter\ifx\csname#1\endcsname\relax
{\bf?}\immediate\write16{#1 ist nicht definiert!}\else\csname#1\endcsname\fi}
\def\expandwrite#1#2{\edef\next{\write#1{#2}}\next}
\def\neverexpand{\noexpand\noexpand\noexpand}
\def\strip#1\ {}
\def\ncite#1{\expandafter\ifx\csname#1\endcsname\relax
{\bf?}\immediate\write16{#1 ist nicht definiert!}\else
\expandafter\expandafter\expandafter\strip\csname#1\endcsname\fi}
\newwrite\AUX
\immediate\openout\AUX=\jobname.aux
\font\eightrm=cmr8\font\sixrm=cmr6
\font\eighti=cmmi8
\font\eightit=cmti8
\font\eightbf=cmbx8
\font\eightcsc=cmcsc10 scaled 833
\def\eightpoint{%
\textfont0=\eightrm\scriptfont0=\sixrm\def\rm{\fam0\eightrm}%
\textfont1=\eighti
\textfont\bffam=\eightbf\def\bf{\fam\bffam\eightbf}%
\textfont\itfam=\eightit\def\it{\fam\itfam\eightit}%
\def\csc{\eightcsc}%
\setbox\strutbox=\hbox{\vrule height7pt depth2pt width0pt}%
\normalbaselineskip=0,8\normalbaselineskip\normalbaselines\rm}
\newcount\absFnno\absFnno1
\write\AUX{\relFnno1}
\newif\ifMARKE\MARKEtrue
{\catcode`\@=11
\gdef\footnote{\ifMARKE\edef\@sf{\spacefactor\the\spacefactor}\/%
$^{\cite{Fn\the\absFnno}}$\@sf\fi
\MARKEtrue
\insert\footins\bgroup\eightpoint
\interlinepenalty100\let\par=\endgraf
\leftskip=0pt\rightskip=0pt
\splittopskip=10pt plus 1pt minus 1pt \floatingpenalty=20000\smallskip
\item{$^{\cite{Fn\the\absFnno}}$}%
\expandwrite\AUX{\neverexpand\ref{Fn\the\absFnno}{\neverexpand\the\relFnno}}%
\global\advance\absFnno1\write\AUX{\advance\relFnno1}%
\bgroup\strut\aftergroup\@foot\let\next}}
\skip\footins=12pt plus 2pt minus 4pt
\dimen\footins=30pc
\output={\plainoutput\immediate\write\AUX{\relFnno1}}
\newcount\Abschnitt\Abschnitt0
\def\beginsection#1. #2 \par{\advance\Abschnitt1%
\vskip0pt plus.10\vsize\penalty-250
\vskip0pt plus-.10\vsize\bigskip\vskip\parskip
\edef\TEST{\number\Abschnitt}
\expandafter\ifx\csname#1\endcsname\TEST\relax\else
\immediate\write16{#1 hat sich geaendert!}\fi
\expandwrite\AUX{\neverexpand\ref{#1}{\TEST}}
\leftline{\marginnote{#1}\bf\number\Abschnitt. \ignorespaces#2}%
\nobreak\smallskip\noindent\SATZ1\GNo0}
\def\Proof:{\par\noindent{\it Proof:}}
\def\Remark:{\ifdim\lastskip<\medskipamount\removelastskip\medskip\fi
\noindent{\bf Remark:}}
\def\Remarks:{\ifdim\lastskip<\medskipamount\removelastskip\medskip\fi
\noindent{\bf Remarks:}}
\def\Definition:{\ifdim\lastskip<\medskipamount\removelastskip\medskip\fi
\noindent{\bf Definition:}}
\def\Example:{\ifdim\lastskip<\medskipamount\removelastskip\medskip\fi
\noindent{\bf Example:}}
\def\Examples:{\ifdim\lastskip<\medskipamount\removelastskip\medskip\fi
\noindent{\bf Examples:}}
\newif\ifmarginalnotes\marginalnotesfalse
\newif\ifmarginalwarnings\marginalwarningstrue

\def\marginnote#1{\ifmarginalnotes\hbox to 0pt{\eightpoint\hss #1\ }\fi}

\def\strutdepth{\dp\strutbox}
\def\Randbem#1#2{\ifmarginalwarnings
{#1}\strut
\setbox0=\vtop{\eightpoint
\rightskip=0pt plus 6mm\hfuzz=3pt\hsize=16mm\noindent\leavevmode#2}%
\vadjust{\kern-\strutdepth
\vtop to \strutdepth{\kern-\ht0
\hbox to \hsize{\kern-16mm\kern-6pt\box0\kern6pt\hfill}\vss}}\fi}

\def\Zitat!{\Randbem{\bf?}{\bf Zitat}}

\newcount\SATZ\SATZ1
\def\proclaim #1. #2\par{\ifdim\lastskip<\medskipamount\removelastskip
\medskip\fi
\noindent{\bf#1.\ }{\it#2}\Par
\ifdim\lastskip<\medskipamount\removelastskip\goodbreak\medskip\fi}
\def\Aussage#1{%
\expandafter\def\csname#1\endcsname##1.{\ifx?##1?\relax\else
\edef\TEST{#1\penalty10000\ \number\Abschnitt.\number\SATZ}
\expandafter\ifx\csname##1\endcsname\TEST\relax\else
\immediate\write16{##1 hat sich geaendert!}\fi
\expandwrite\AUX{\neverexpand\ref{##1}{\TEST}}\fi
\proclaim {\marginnote{##1}\number\Abschnitt.\number\SATZ. #1\global\advance\SATZ1}.}}
\Aussage{Theorem}
\Aussage{Proposition}
\Aussage{Corollary}
\Aussage{Lemma}
\font\la=lasy10
\def\strich{\hbox{$\vcenter{\hbox
to 1pt{\leaders\hrule height -0,2pt depth 0,6pt\hfil}}$}}
\def\dashedrightarrow{\hbox{%
\hbox to 0,5cm{\leaders\hbox to 2pt{\hfil\strich\hfil}\hfil}%
\kern-2pt\hbox{\la\char\string"29}}}

\def\Bindestrich{\penalty10000-\hskip0pt}
\let\_=\Bindestrich
\def\.{{\sfcode`.=1000.}}
\def\Links#1{\llap{$\scriptstyle#1$}}
\def\Rechts#1{\rlap{$\scriptstyle#1$}}
\def\Par{\par}
\def\:={\mathrel{\raise0,9pt\hbox{.}\kern-2,77779pt
\raise3pt\hbox{.}\kern-2,5pt=}}
\def\=:{\mathrel{=\kern-2,5pt\raise0,9pt\hbox{.}\kern-2,77779pt
\raise3pt\hbox{.}}} \def\mod{/\mskip-5mu/}
\def\into{\hookrightarrow}
\def\pfeil{\rightarrow}
\def\untenPf{\downarrow}
\def\Pfeil{\longrightarrow}
\def\pf#1{\buildrel#1\over\rightarrow}
\def\Pf#1{\buildrel#1\over\longrightarrow}

\def\Ugleich{\hbox{$\cup$\kern.5pt\vrule depth -0.5pt}}
\def\|#1|{\mathop{\rm#1}\nolimits}
\def\<{\langle}
\def\>{\rangle}
\let\Times=\times
\def\times{\mathop{\Times}}
\let\Otimes=\otimes
\def\otimes{\mathop{\Otimes}}
\catcode`\@=11
\def\hex#1{\ifcase#1 0\or1\or2\or3\or4\or5\or6\or7\or8\or9\or A\or B\or
C\or D\or E\or F\else\message{Warnung: Setze hex#1=0}0\fi}
\def\fontdef#1:#2,#3,#4.{%
\alloc@8\fam\chardef\sixt@@n\FAM
\ifx!#2!\else\expandafter\font\csname text#1\endcsname=#2
\textfont\the\FAM=\csname text#1\endcsname\fi
\ifx!#3!\else\expandafter\font\csname script#1\endcsname=#3
\scriptfont\the\FAM=\csname script#1\endcsname\fi
\ifx!#4!\else\expandafter\font\csname scriptscript#1\endcsname=#4
\scriptscriptfont\the\FAM=\csname scriptscript#1\endcsname\fi
\expandafter\edef\csname #1\endcsname{\fam\the\FAM\csname text#1\endcsname}
\expandafter\edef\csname hex#1fam\endcsname{\hex\FAM}}
\catcode`\@=12 

\fontdef Ss:cmss10,,.
\fontdef Fr:eufm10,eufm7,eufm5.
\def\fa{{\Fr a}}
\def\fb{{\Fr b}}

\def\fg{{\Fr g}}
\def\fh{{\Fr h}}

\def\fk{{\Fr k}}
\def\fl{{\Fr l}}
\def\fm{{\Fr m}}

\def\fp{{\Fr p}}

\def\ft{{\Fr t}}
\def\fu{{\Fr u}}

\fontdef bbb:msbm10,msbm7,msbm5.
\fontdef mbf:cmmib10,cmmib7,.
\fontdef msa:msam10,msam7,msam5.
\def\CC{{\bbb C}}

\def\cA{{\cal A}}\def\cB{{\cal B}}
\def\cG{{\cal G}}\def\cH{{\cal H}}
\def\cL{{\cal L}}
\def\cM{{\cal M}}

\def\cU{{\cal U}}

\mathchardef\leer=\string"0\hexbbbfam3F
\mathchardef\subsetneq=\string"3\hexbbbfam24
\mathchardef\semidir=\string"2\hexbbbfam6E
\mathchardef\dirsemi=\string"2\hexbbbfam6F
\mathchardef\haken=\string"2\hexmsafam78
\mathchardef\auf=\string"3\hexmsafam10
\let\OL=\overline
\def\overline#1{{\hskip1pt\OL{\hskip-1pt#1\hskip-.3pt}\hskip.3pt}}

\def\Cq{{\overline{C}}}

\def\Fq{{\overline{F}}}
\def\Gq{{\overline{G}}}

\def\Xq{{\overline{X}}}

\def\Zq{{\overline{Z}}}
%
\newdimen\Parindent
\Parindent=\parindent


\abovedisplayskip 9.0pt plus 3.0pt minus 3.0pt
\belowdisplayskip 9.0pt plus 3.0pt minus 3.0pt
\newdimen\Grenze\Grenze2\Parindent\advance\Grenze1em
\newdimen\Breite
\newbox\DpBox
\def\NewDisplay#1$${\Breite\hsize\advance\Breite-\hangindent
\setbox\DpBox=\hbox{\hskip2\Parindent$\displaystyle{#1}$}%
\ifnum\predisplaysize<\Grenze\abovedisplayskip\abovedisplayshortskip
\belowdisplayskip\belowdisplayshortskip\fi
\global\futurelet\nexttok\WEITER}
\def\WEITER{\ifx\nexttok\qed\expandafter\leftQEDdisplay
\else\leftdisplay\fi}
\def\leftdisplay{\hskip-\hangindent\leftline{\box\DpBox}$$}
\def\leftQEDdisplay{\hskip-\hangindent
\line{\copy\DpBox\hfill\lower\dp\DpBox\copy\QEDbox}%
\belowdisplayskip0pt$$\bigskip\let\nexttok=}
\everydisplay{\NewDisplay}
\newcount\GNo\GNo=0
\newcount\maxEqNo\maxEqNo=0
\def\eqno#1{
\global\advance\GNo1
\edef\FTEST{$(\number\Abschnitt.\number\GNo)$}
\ifx?#1?\relax\else
\ifnum#1>\maxEqNo\global\maxEqNo=#1\fi%
\expandafter\ifx\csname E#1\endcsname\FTEST\relax\else
\immediate\write16{E#1 hat sich geaendert!}\fi
\expandwrite\AUX{\neverexpand\ref{E#1}{\FTEST}}\fi
\llap{\hbox to 40pt{\marginnote{#1}\FTEST\hfill}}}

\catcode`@=11
\def\eqalignno#1{\null\vcenter{\openup\jot\m@th\ialign{\eqno{##}\hfil
&\strut\hfil$\displaystyle{##}$&$\displaystyle{{}##}$\hfil\crcr#1\crcr}}\,}
\catcode`@=12

\newbox\QEDbox
\newbox\nichts\setbox\nichts=\vbox{}\wd\nichts=2mm\ht\nichts=2mm
\setbox\QEDbox=\hbox{\vrule\vbox{\hrule\copy\nichts\hrule}\vrule}
\def\qed{\leavevmode\unskip\hfil\null\nobreak\hfill\copy\QEDbox\medbreak}
\newdimen\HIindent
\newbox\HIbox
\def\setHI#1{\setbox\HIbox=\hbox{#1}\HIindent=\wd\HIbox}
\def\HI#1{\par\hangindent\HIindent\hangafter=0\noindent\leavevmode
\llap{\hbox to\HIindent{#1\hfil}}\ignorespaces}

\newdimen\maxSpalbr
\newdimen\altSpalbr
\newcount\Zaehler


\newif\ifxxx

{\catcode`/=\active

\gdef\beginrefs{%
\xxxfalse
\catcode`/=\active
\def/{\string/\ifxxx\hskip0pt\fi}
\def\TText##1{{\xxxtrue\tt##1}}
\expandafter\ifx\csname Spaltenbreite\endcsname\relax
\def\Spaltenbreite{1cm}\immediate\write16{Spaltenbreite undefiniert!}\fi
\expandafter\altSpalbr\Spaltenbreite
\maxSpalbr0pt
\gdef\alt{}
\def\\##1\relax{%
\gdef\neu{##1}\ifx\alt\neu\global\advance\Zaehler1\else
\xdef\alt{\neu}\global\Zaehler=1\fi\xdef\SigText{##1\the\Zaehler}}
\def\L|Abk:##1|Sig:##2|Au:##3|Tit:##4|Zs:##5|Bd:##6|S:##7|J:##8|xxx:##9||{%
\def\SigText{##2}\global\setbox0=\hbox{##2\relax}
\edef\TEST{[\SigText]}
\expandafter\ifx\csname##1\endcsname\TEST\relax\else
\immediate\write16{##1 hat sich geaendert!}\fi
\expandwrite\AUX{\neverexpand\ref{##1}{\TEST}}
\setHI{[\SigText]\ }
\ifnum\HIindent>\maxSpalbr\maxSpalbr\HIindent\fi
\ifnum\HIindent<\altSpalbr\HIindent\altSpalbr\fi
\HI{\marginnote{##1}[\SigText]}
\ifx-##3\relax\else{##3}: \fi
\ifx-##4\relax\else{##4}{\sfcode`.=3000.} \fi
\ifx-##5\relax\else{\it ##5\/} \fi
\ifx-##6\relax\else{\bf ##6} \fi
\ifx-##8\relax\else({##8})\fi
\ifx-##7\relax\else, {##7}\fi
\ifx-##9\relax\else, \TText{##9}\fi\Par}
\def\B|Abk:##1|Sig:##2|Au:##3|Tit:##4|Reihe:##5|Verlag:##6|Ort:##7|J:##8|xxx:##9||{%
\def\SigText{##2}\global\setbox0=\hbox{##2\relax}
\edef\TEST{[\SigText]}
\expandafter\ifx\csname##1\endcsname\TEST\relax\else
\immediate\write16{##1 hat sich geaendert!}\fi
\expandwrite\AUX{\neverexpand\ref{##1}{\TEST}}
\setHI{[\SigText]\ }
\ifnum\HIindent>\maxSpalbr\maxSpalbr\HIindent\fi
\ifnum\HIindent<\altSpalbr\HIindent\altSpalbr\fi
\HI{\marginnote{##1}[\SigText]}
\ifx-##3\relax\else{##3}: \fi
\ifx-##4\relax\else{##4}{\sfcode`.=3000.} \fi
\ifx-##5\relax\else{(##5)} \fi
\ifx-##7\relax\else{##7:} \fi
\ifx-##6\relax\else{##6}\fi
\ifx-##8\relax\else{ ##8}\fi
\ifx-##9\relax\else, \TText{##9}\fi\Par}
\def\Pr|Abk:##1|Sig:##2|Au:##3|Artikel:##4|Titel:##5|Hgr:##6|Reihe:{%
\def\SigText{##2}\global\setbox0=\hbox{##2\relax}
\edef\TEST{[\SigText]}
\expandafter\ifx\csname##1\endcsname\TEST\relax\else
\immediate\write16{##1 hat sich geaendert!}\fi
\expandwrite\AUX{\neverexpand\ref{##1}{\TEST}}
\setHI{[\SigText]\ }
\ifnum\HIindent>\maxSpalbr\maxSpalbr\HIindent\fi
\ifnum\HIindent<\altSpalbr\HIindent\altSpalbr\fi
\HI{\marginnote{##1}[\SigText]}
\ifx-##3\relax\else{##3}: \fi
\ifx-##4\relax\else{##4}{\sfcode`.=3000.} \fi
\ifx-##5\relax\else{In: \it ##5}. \fi
\ifx-##6\relax\else{(##6)} \fi\PrII}
\def\PrII##1|Bd:##2|Verlag:##3|Ort:##4|S:##5|J:##6|xxx:##7||{%
\ifx-##1\relax\else{##1} \fi
\ifx-##2\relax\else{\bf ##2}, \fi
\ifx-##4\relax\else{##4:} \fi
\ifx-##3\relax\else{##3} \fi
\ifx-##6\relax\else{##6}\fi
\ifx-##5\relax\else{, ##5}\fi
\ifx-##7\relax\else, \TText{##7}\fi\Par}
\bgroup
\baselineskip12pt
\parskip2.5pt plus 1pt
\hyphenation{Hei-del-berg}
\sfcode`.=1000
\beginsection References. References

}}

\def\endrefs{%
\expandwrite\AUX{\neverexpand\ref{Spaltenbreite}{\the\maxSpalbr}}
\ifnum\maxSpalbr=\altSpalbr\relax\else
\immediate\write16{Spaltenbreite hat sich geaendert!}\fi
\egroup\write16{Letzte Gleichung: E\the\maxEqNo}}



\Aussage{Example}
\def\r#1{{\it\romannumeral#1)}}
\def\ol{\overline}
\def\punkt{{\hbox{\bf.}}}
\def\dotunion{\mathop{\textstyle\buildrel\punkt\over\bigcup}}

\title{A connectedness property of algebraic moment maps}

\removelastskip\bigskip \vskip \parskip
\leftline{\it For Claudio Procesi}

\beginsection Introduction. Introduction

Let $K$ be a connected compact Lie group and $M$ a Hamiltonian
$K$\_manifold, i.e., a symplectic $K$\_manifold equipped with a moment
map $\mu:M\pfeil\fk^*:=(\|Lie|K)^*$. A theorem of Kirwan (implicitly
in \cite{Kir}) asserts: if $M$ is connected and compact then the level
sets of $\mu$ are connected. The purpose of this note is to prove such
a statement in the category of algebraic varieties.

First, we reformulate Kirwan's theorem: consider the map
$\psi:M\pfeil\fk^*/K$ which is the composition of $\mu$ with the
quotient map. For a point $x\in\fk^*$ let $H=K_x$ be its isotropy group and
$y=Kx\in\fk^*/K$ its orbit. Then the fiber $\psi^{-1}(y)$ is
isomorphic to the fiber product $K\times^H\mu^{-1}(x)$. Thus, since
both $K/H$ and $H$ are connected, the connectedness of the fibers of
$\mu$ is equivalent to the connectedness of the fibers of $\psi$. This
formulation is more suitable for the algebraic category.

Let $G$ be a connected reductive group (everything over $\CC$) and $Z$
a Hamiltonian $G$\_variety with moment map
$\mu:Z\pfeil\fg^*=(\|Lie|G)^*$. Let $\fg^*\mod G:=\|Spec|\CC[\fg^*]^G$
be the categorical quotient and let $\tilde\psi:Z\pfeil\fg^*\mod G$
be the composition of $\mu$ with the quotient map. Since the latter is
not an orbit map the connection between fibers of $\mu$ and
$\tilde\psi$ is much more loose than in the differential category and
we concentrate on $\tilde\psi$ from now on.

The morphism $\tilde\psi$ is still not the right map, since sometimes
not even its generic fibers are connected. An example is the action of
$G=SL_2(\CC)$ on $Z=\CC^2\times(\CC^2)^*$ (this is the cotangent
bundle of $\CC^2$). Then $\fg^*\mod G=\A^1$ and
$\tilde\psi(u,\alpha)=\alpha(u)^2$. In particular, the generic fiber
of $\tilde\psi$ has two connected components. This is remedied by
looking at the map $\psi(u,\alpha):=\alpha(u)$ instead. Then
$\tilde\psi$ is the composition of $\psi$ with the finite map
$\A^1\pfeil\A^1:z\mapsto z^2$, the latter being responsible for the
disconnected fibers.

This construction can be generalized as follows: the morphism
$\tilde\psi$ induces an homomorphism of algebras
$$\eqno{}
\tilde\psi^*:\ \CC[\fg^*\mod G]\pfeil\CC[Z].
$$
Let $R\subseteq\CC[Z]$ be the integral closure of the image and
$L:=\|Spec|R$. Then $\tilde\psi$ factors through $\psi:Z\pfeil L$. By
construction, the generic fibers of $\psi$ are now connected, even
irreducible. Furthermore, we expect that the behavior of $\psi$ is
dramatically better than that. For example, when $Z$ is connected and
affine, then it is hoped that $L$ is smooth and $\psi$ is faithfully
flat with reduced, connected fibers. If that were true then most bad
behavior of $\tilde\psi$ is on account of the finite morphism
$L\pfeil\fg^*\mod G$.

So far, the theory of algebraic Hamiltonian varieties is not developed
enough to convert this hope into proofs, but for one class we know a
lot. That is when $Z$ is the cotangent bundle of a smooth $G$\_variety
$X$. Then one can prove, \cite{WuM}, that $L$ is not only smooth but
even an affine space. More precisely, $L=\fa^*/W_X$ where $\fa^*$ is a
finite dimensional vector space and $W_X$ a finite reflection
group. In this setting, the main result of this paper is:

\Theorem ConnThm. Let $X$ be a connected, smooth $G$\_variety and
$Z:=T^*_X$ its cotangent bundle. Then all fibers of
$\psi:Z\pfeil\fa^*/W_X$ are connected.

Apart from its intrinsic interest this theorem is useful for further
investigation of the fibers of $\psi$. For example, in \cite{WHM} it
is used to prove that most fibers of $\psi$ are reduced, provided $X$
is affine. In turn, that latter result is crucial in the investigation
of so\_called collective functions on Hamiltonian manifolds in the
differentiable (!) category (also in \cite{WHM}).

The strategy of the proof of \cite{ConnThm} is as follows: the core
argument is an application of Zariski's connectedness theorem, a
theorem which works only for proper morphisms. For that reason we
construct a partial compactification of $Z$. Since this works only
when $X$ is homogeneous the proof splits in two parts: 1. proof of the
connectedness theorem for homogeneous $X$ and 2. reduction of the
general case to the homogeneous case.

In both parts another key ingredient are certain irreducible
subvarieties of the fibers of $\psi$ which are used to ``tie'' its
different parts together, thereby proving connectivity (cf. \cite{BBB}
and the construction of $E$ in \S\cite{generalcase}).

\Remark: This paper is an expanded version of \cite{WHM} \S6.
\medskip

\noindent{\bf Notation:} All varieties are defined over an
algebraically closed field of characteristic zero, which is denoted by
$\CC$. Throughout the paper, $G$ is a connected reductive group. The
Lie algebras of the algebraic groups $G$, $B$, $U$, $H$,~$\ldots$ are
denoted by the corresponding lower case fraktur letters $\fg$, $\fb$,
$\fu$, $\fh$, $\ldots$ For an affine $G$\_variety $X$ let $X\mod
G:=\|Spec|\CC[X]^G$ be the categorical quotient. For a subspace
$U\subseteq V$ let $U^\perp\subseteq V^*$ denote its annihilator.

\beginsection Review. The moment map on cotangent bundles

In this section we review some results of  \cite{WuM} about the geometry
of the moment map on a cotangent bundle.

Let $X$ be a smooth connected $G$\_manifold. Let $B\subseteq G$ be a
Borel subgroup with maximal unipotent subgroup $U$ and maximal torus
$T$. Then one can define the following important numerical invariants
of $X$:
$$\eqno{}
n:=\|dim|X,\quad n_u:=\|max|\limits_{x\in X}\ \|dim|Ux,\quad
n_b:=\|max|\limits_{x\in X}\ \|dim|Bx.
$$
The difference $c:=n-n_b$ is called the {\it complexity} of $X$ while
$r:=n_b-n_u$ is its {\it rank}. Consider
$$\eqno{}
X_0:=\{x\in X\mid \|dim| Ux=n_u,\|dim|Bx=n_b\}.
$$
This is a dense open $B$\_stable subset of $X$. We define the
following subbundles of the cotangent bundle $T^*_{X_0}$:
$$\eqno{}
\cB:=\dotunion_{x\in X_0}(\fb x)^\perp\subseteq\cU:=\dotunion_{x\in X_0}(\fu
x)^\perp.
$$
The fibers of the quotient are
$$\eqno{12}
(\cU/\cB)_x=\cU_x/\cB_x=(\fu x)^\perp/(\fb x)^\perp=(\fb x/\fu
x)^*=(\fb/(\fu+\fb_x))^*\subseteq(\fb/\fu)^*.
$$
Therefore we can identify this fiber with a subspace $\fa_x^*$ of
the ``abstract'' dual Cartan subalgebra $\ft^*=(\fb/\fu)^*$ . By
construction, the ranks of the vector bundles $\cB$ and $\cU$ are $c$
and $c+r$, respectively. Hence we have $\|dim|\fa^*_x=r$.

\Lemma independent. {\rm(\cite{WuM}~6.3)} The space $\fa^*=\fa^*_x$
depends neither on the choice of $x\in X_0$ nor on the choice of a
Borel subgroup $B$. In particular, we get a morphism $\cU\auf\fa^*$
which induces a trivialization $\cU/\cB\pf\sim\fa^*\times X_0$.

The cotangent bundle $T^*_X$ is a Hamiltonian $G$\_variety, i.e., it
carries a natural $G$\_invariant symplectic structure and is equipped
with a moment map.  More precisely, let $\pi:T^*_X\pfeil X$ be the
natural projection.  Then the moment map is
$$\eqno{}
\mu:T^*_X\pfeil \fg^*:\alpha\mapsto[\xi\mapsto\alpha(\xi\pi(\alpha))].
$$
Actually, we are more interested in the composition
$\tilde\psi:T^*_X\pf\mu\fg^*\pfeil\fg^*\mod G$. It induces a
homomorphisms of algebras
$$\eqno{}
{\tilde\psi}^*:\CC[\fg^*\mod G]\pfeil\CC[T^*_X].
$$
Let $L_X$ be the spectrum of the integral closure of the
$\tilde\psi^*\CC[\fg^*]$ in $\CC[T^*_X]$. Then we obtain a
factorization of $\tilde\psi$
$$\eqno{}
T^*_X\pf\psi L_X\pf\eta\fg^*\mod G
$$
where $\eta$ is a finite morphism and $\psi$ is dominant. Moreover,
$\psi$ is universal with these properties. An important property of
$\psi$ is that its generic fibers are irreducible.

\Theorem restriction. {\rm(\cite{WuM}~6.2,~6.6a,b)} The restriction of
$\psi$ to $\cU$ factors through $\fa^*$:
$$\eqno{11}
\cU\pfeil\fa^*\pf\pi L_X.
$$
The morphism $\pi$ is finite and surjective. More precisely, there is
a finite reflection group $W_X\subseteq GL(\fa^*)$ (the little Weyl
group of $X$) such that $\pi$ induces an isomorphism $\fa^*/W_X\pf\sim
L_X$. In particular, $L_X$ is an affine space.

\noindent From now on, we identify $L_X$ with $\fa^*/W_X$. It will
also be convenient to use a less canonical way of stating
\cite{E11}. Let $x\in X_0$ and let $\fa'$ be a complement of $(\fb
x)^\perp\subseteq T^*_{X,x}$ in $(\fu x)^\perp\subseteq
T^*_{X,x}$. Then, by \cite{E12} we can identify $\fa'$ with $\fa^*$
and we have a commutative diagram
$$\eqno{}
\matrix{\fa'&\into&T^*_X\cr
\untenPf\Rechts\sim&&\untenPf\Rechts\psi\cr
\fa^*&\auf&\fa^*/W_X\cr}
$$
In other words, $\fa'$ is ``almost'' a section of $\psi$.

\Theorem faithflat. {\rm(\cite{WuM}~6.6c)} The morphism
$\psi:T^*_X\pfeil\fa^*/W_X$ is faithfully flat. In particular, all
fibers of $\psi$ are of pure codimension $r$.

Now we specialize everything to the case where $X=G/H$ is a
homogeneous space. Then we can write $T^*_X=G\times^H\fh^\perp$ and all
information about $\psi$ is contained in the restriction
$$\eqno{}
\psi_\fh:\fh^\perp\pfeil\fa^*/W_X.
$$

\Theorem equilie. Let $X=G/H$ be a homogeneous variety.
\item{\r1} The rank $r$, the complexity $c$, and the little Weyl group
$W_{G/H}$ depend only on the Lie algebra $\fh$.\Par
\item{\r2} The morphism $\psi_\fh$ is faithfully flat. All of its
fibers are pure of codimension $r$. The generic fibers are
irreducible.\Par
\item{\r3} Let $B$ be a Borel subgroup of $G$ and $g:=\|dim|G$. Then
$$\eqno{6}
\|dim|(\fb+\fh)\le g-c\quad\hbox{and}\quad\|dim|(\fu+\fh)\le g-c-r.
$$
Moreover, equality holds for an open set of Borel subgroups.\Par
\item{\r4} Assume that equality holds in \cite{E6}. Let $\fa'$ be any
complement of $\fh^\perp\cap\fb^\perp$ in
$\fh^\perp\cap\fu^\perp$. Then one can identify $\fa'$ with $\fa^*$
and the restriction of $\psi_\fh$ to $\fa'$ is the quotient map
$\fa^*\pfeil\fa^*/W_X$.\Par

\Proof: \r1 The connected component $H^0$ of $H$ depends only on
$\fh$. Moreover, the morphism $X^0:=G/H^0\pfeil X=G/H$ is
\'etale. This implies that dimension, complexity, and rank of $X^0$
and $X$ are the same (a direct consequence of the definitions). The same
holds for the little Weyl group (see \cite{WuM}~6.5.3).

\r2 By \cite{faithflat}, the morphism
$\psi_0:T^*_{X^0}\pfeil\fa^*/W_{X^0}=\fa^*/W_X$ is surjective and has
fibers of codimension $r$. Moreover, the generic fibers are
irreducible. Since, for $u\in\fa^*/W_X$ we have
$\psi_0^{-1}(u)=G\times^{H^0}\psi_\fh^{-1}(u)$, the same holds for
$\psi_\fh$. It follows that $\psi_\fh$ is faithfully flat (\cite{EGA}
\S15.4.2).

\r3 Let $x_0=eH\in G/H$ be the base point. Then
$\|dim|(\fb+\fh)=\|dim|Bx_0+\|dim|\fh$ and
$\|dim|(\fu+\fh)=\|dim|Ux_0+\|dim|\fh$. Thus, \cite{E6} follows from
the definition of $c$ and $r$. Moreover, there is a point $gx_0\in
G/H$ such that $Bgx_0$ and $Ugx_0$ have maximal dimension. Thus, if we
replace $B$ by $g^{-1}Bg$ we get equalities in \cite{E6}.

\r4 We have
$$\eqno{}
\fa'\pf\sim(\fh^\perp\cap\fu^\perp)/(\fh^\perp\cap\fb^\perp)=
(\fh+\fu)^\perp/(\fh+\fb)^\perp=[(\fh+\fb)/(\fh+\fu)]^*=[\fb x_0/\fu
x_0]^*.
$$
Now the assertion follows from \cite{restriction}.\qed

\beginsection construction. A partial compactification: construction

Ultimately, we would like to invoke Zariski's connectedness theorem
which allows one to deduce the connectedness of all fibers from that
of the generic fibers. The essential prerequisite for this theorem is
that the morphism is proper. Therefore, we construct first a partial
compactification of the cotangent bundle which renders the moment map
proper. For this we need that $X=G/H$ is
homogeneous. So we assume this until further notice.

Consider the Lie algebra $\fh$ as a point of the Grassmannian
$\|Gr|(\fg)$ of all subspaces of $\fg$. Let $Y\subseteq\|Gr|(\fg)$ be
the closure of the orbit $G\cdot\fh$. Recall that a subalgebra of
$\fg$ is called {\it algebraic\/} if it is the Lie algebra of a closed
subgroup of $G$. By \cite{equilie}\r1 we may speak of the rank, the
complexity, and the little Weyl group of an algebraic Lie subalgebra
of $\fg$.

\Lemma alglie. Every point $\fm\in Y$ represents an algebraic
subalgebra of $\fg$ which has the same dimension, complexity, and rank
as $\fh$.

\Proof: First, we claim that there is a smooth affine curve $C$, a
point $c_0\in C$, and morphisms $\alpha':C':=C\setminus\{c_0\}\pfeil
G$, $\alpha:C\pfeil Y$ such that 
$$\eqno{2}
\alpha(c)=\cases{\alpha'(c)\cdot\fh&for $c\ne c_0$;\cr
\fm&for $c=c_0$.\cr}
$$
In fact, let $G\into G'$ be any completion and $\Gq\subseteq G'\times
Y$ the closure of the set $\{(g,g\cdot\fh)\mid g\in G\}$. Then
$$\eqno{3}
G\into\Gq:g\mapsto(g,g\cdot\fh)
$$
is again a completion of $G$. Let $\gamma:\Gq\pfeil Y$ be the second
projection. This morphism is dominant and proper, hence
surjective. Choose $x\in\Gq$ with $\gamma(x)=\fm$. Since $G$ is open
and dense in $\Gq$ there is an affine curve $C_1\subseteq\Gq$ with
$x\in C_1$ and $C_1\setminus\{x\}\subseteq G$. Let
$\CS\pfeil C_1$ be the normalization of $C$ and $S$ the preimage of
$x$. This is a non\_empty finite subset of $\CS$. The desired curve
$C$ is obtained by removing from $\CS$ all points but one of $S$. The
remaining point of $S$ is called $c_0$. The morphism $\alpha$ is just
the composition $C\pfeil\Gq\pfeil Y$. Moreover, the image of $C'$ in
$\Gq$ is contained in $G$ giving $\alpha':C'\pfeil G$. Finally,
\cite{E2} follows from \cite{E3} which proves the claim.

Using this curve $C$, consider the trivial group schemes $\cG:=G\times
C/C$ and $\cG':=G\times
C'/C'$. Then $\cG'$ contains the subgroup scheme
$$\eqno{}
\cH':=\dotunion_{c\in C'}\alpha'(c)H\alpha'(c)^{-1}=
\{(\alpha'(c)\,h\,\alpha'(c)^{-1},c)\mid h\in H,c\in C'\}\cong H\times C'.
$$
Let $\cH$ be the closure of $\cH'$ in $\cG$. Then every irreducible
component of $\cH$ maps dominantly to $C$ which implies that $\cH$ is
flat over $C$; see e.g. \cite{Har} 9.7. Therefore, the same holds true
for the fiber product $\cH^2:=\cH\times_C\cH$ which implies that
$\cH'\times_C\cH'$ is dense in $\cH^2$. In turn, the multiplication
map $\cG\times_C \cG\pfeil\cG$ maps $\cH^2$ into $\cH$. In other
words, $\cH$ is a flat subgroup scheme of $\cG$. In particular, the
fiber $H_c:=\cH\times_C\{c\}$ is a subgroup of $G$ for every $c\in
C$. Any flat group scheme is smooth (Cartier, see
\cite{DG}). Therefore, the Lie algebras $\fh_c:=\|Lie|H_c$ form a
vector bundle over $C$. By construction we have
$\fh_c=\alpha'(c)\cdot\fh$ for all $c\ne c_0$. Thus, by \cite{E2}, we
have $\fh_c=\alpha(c)$ for all $c$. In particular,
$\fh_{c_0}=\fm$. This shows that $\fm=\|Lie|H_{c_0}$ is algebraic.

Finally, it is known that the quotient $\cG/\cH$ exists
(\cite{An}~4.C). This is a flat deformation of homogeneous
$G$\_varieties with fibers $G/H_c$. Thus dimension, complexity and
rank don't change (\cite{WuM}~2.5) which proves the second
assertion.\qed

Now we construct a certain partial compactification of $T^*_X$. First,
choose an equivariant embedding of $X=G/H$ into a projective space:
$X=G/H\into\P^N$. Using the orbit map
$$\eqno{4}
G/H\pfeil\|Gr|(\fg):gH\mapsto g\cdot\fh.
$$
we get a diagonal embedding of $G/H$ in $\P^N\times\|Gr|(\fg)$. Let
$X'$ be its closure. This is a projective, possibly singular
$G$\_variety. Let $\Xq\pfeil X'$ be an equivariant desingularization
(see \cite{AHV}). Then we have a smooth completion $G/H\into\Xq$
such that the orbit map \cite{E4} extends to a morphism
$\Xq\pfeil\|Gr|(\fg)$.

Let $\cM\pfeil\Xq$ be the pull\_back of the tautological vector bundle
over $\|Gr|(\fg)$. For $x\in\Xq$ let $\fm_x:=\cM_x$ be its fiber over
$x$. Then, by construction, we have $\fm_x=g\cdot\fh$ whenever
$x=gH\in G/H$. Thus, by \cite{alglie}, every $\fm_x$ is an
algebraic Lie subalgebra of $\fg$ which has the same dimension,
complexity, and rank as $\fh$.

By construction, the vector bundle $\cM$ is a subbundle of the trivial
vector bundle $\fg\times\Xq$. Let $\Zq\pfeil\Xq$ be its annihilator in
the dual bundle $\fg^*\times\Xq$. This means that each fiber $\Zq_x$
equals $\fm_x^\perp\subseteq\fg^*$. We claim that $\Zq$ is a partial
compactification of $Z:=T^*_X$. In fact, the restriction of $\Zq$ to
the open orbit $G/H$ is $G\times^H\fh^\perp$ which equals the
cotangent bundle over $X$. Note that $\Zq$ is equipped with a natural
$G$\_equivariant morphism
$$\eqno{5}
\ol\mu:\Zq\into\fg^*\times\Xq\auf\fg^*
$$
which, on each fiber, is simply the natural embedding
$\fm_x^\perp\into\fg^*$. This means, that $\ol\mu$ is an
extension of the moment map $\mu:Z\pfeil\fg^*$ to $\Zq$. The point is
now that $\ol\mu$ is a {\it proper\/} morphism as one sees from
the factorization \cite{E5} and the fact that $\Xq$ is complete.

\beginsection properties. A partial compactification: properties

Now, a first step to our goal is the following lemma:

\Lemma EEE. The morphism $\psi:Z\pfeil\fa^*/W_X$ extends to a morphism
$\ol\psi:\Zq\pfeil\fa^*/W_X$. Moreover, all fibers of
$\ol\psi$ are connected.

\Proof: Let $\Gamma_\psi\subseteq Z\times\fa^*/W_X$ be the graph of
$\psi$ and let $\ol\Gamma_\psi$ be its closure in
$\Zq\times\fa^*/W_X$. Then the projection $\ol\Gamma_\psi\pfeil\Zq$ is
an isomorphism over $X$, hence birational. Now observe that the
composition $Z\pfeil\fa^*/W_X\pfeil\fg^*\mod G$ does have an extension
to $\Zq$, namely the composition $\Zq\pfeil\fg^*\pfeil\fg^*\mod
G$. This implies that $\Gamma_\psi$, hence $\ol\Gamma_\psi$, is
contained in the fiber product $\Gamma':=\Zq\times_{\fg^*\mod
G}\fa^*/W_X$. Since $\fa^*/W_X\pfeil\fg^*\mod G$ is a finite morphism,
so are $\Gamma'\pfeil\Zq$ and $\ol\Gamma_\psi\pfeil\Zq$. Thus, the
latter is an isomorphism since $\Zq$ is normal (even smooth). This
shows that $\ol\psi$ exists.

To show connectedness of fibers, we need another construction of
\cite{WuM}. Consider the moment map $\mu:Z\pfeil\fg^*$ and let $M$ be
the spectrum of the integral closure of $\mu^*\CC[\fg^*]$ in
$\CC[Z]$. Then $\mu$ factors as $Z\pf\phi M\pfeil\fg^*$. As above, one
shows that $\phi$ extends to a morphism $\ol\phi:\Zq\pfeil M$ which
factors $\ol\mu$. This implies that also $\ol\phi$ is proper. By
construction of $M$, the generic fibers of $\phi$ are
irreducible. Hence the same holds for $\ol\phi$. Since $M$ is normal
we can apply Zariski's connectedness theorem (\cite{Har}~11.3) and
conclude that all fibers of $\ol\phi$ are connected. On the other
hand, we have a commutative diagram
$$\eqno{}
\matrix{\Zq&\pf{\ol\phi}&M\cr
&\Links{\ol\psi}\searrow&\downarrow\Rechts{\pi}\cr
&&\fa^*/W_X}
$$
If $f$ is a $G$\_invariant regular function on $Z$ which is integral
over $\CC[\fg^*]$ then it is also integral over $\CC[\fg^*]^G$ (just
apply the Reynolds operator to an integral equation of $f$). This
means that $\pi$ is just the categorical quotient by $G$. Thus every
fiber of $\pi$ contains a unique (connected) closed orbit which
implies that $\pi$ has connected fibers, as well. Choose
$u\in\in\fa^*/W_X$. By applying the following lemma to the morphism
$\ol\psi^{-1}(u)\pfeil\pi^{-1}(u)$ we conclude that $\ol\psi^{-1}(u)$
is connected.\qed

\Lemma openclosed. Let $\phi:X\pfeil Y$ be a surjective morphism 
whose fibers are all connected. Assume $\phi$ is either closed
(e.g., proper) or open (e.g., flat). Assume moreover that $Y$ is
connected. Then $X$ is connected as well.

\Proof: Assume that $\phi$ is closed/open. Suppose that $X$ is not
connected. Then $X=X_1\dot\cup X_2$ is a disjoint union of non\_empty
closed/open subsets. Thus, $Y=\phi(X_1)\cup\phi(X_2)$
is a union of non\_empty closed/open subsets. Since $Y$ is connected, this
union cannot be disjoint. Choose $y\in\phi(X_1)\cap\phi(X_2)$. Then
the fiber $F=\phi^{-1}(y)$ is the union of two disjoint non\_empty
closed/open subsets, namely $F\cap X_1$ and $F\cap X_2$, which contradicts
the connectedness of $F$.
\qed

Next, we investigate the restriction of $\ol\psi$ to the fibers
$\Zq_x=\fm_x^\perp$. Let $c$ and $r$ be the complexity and rank of
$G/H$. For a fixed Borel subgroup $B$ let $\Xq_0$ be the set of
$x\in\Xq$ such that
$$\eqno{}
\|dim|(\fb+\fm_x)=g-c\quad\hbox{and}\quad\|dim|(\fu+\fm_x)=g-c-r.
$$
Since $\fm_x$ has also complexity $c$ and rank $r$ (\cite{alglie})
this is an open subset of $\Xq$ which intersects each $G$\_orbit
non\_trivially. Let
$$\eqno{}
\ol\cB=\dotunion_{x\in\Xq_0}(\fb+\fm_x)^\perp
\quad\hbox{and}\quad
\ol\cU=\dotunion_{x\in\Xq_0}(\fu+\fm_x)^\perp.
$$
Then $\ol\cB\subseteq\ol\cU\subseteq\Zq$ are sub\_vector
bundles.  Consider the fiber of the quotient
$\fa^*_x=\ol\cU_x/\ol\cB_x= (\fb/\fu+\fb\cap\fm_x)^*$ of
$(\fb/\fu)^*=\ft^*$. For $x\in X\cap\Xq_0$ it is always the same space
$\fa^*$ (\cite{independent}). Thus, by continuity, $\fa^*_x=\fa^*$ for
all $x\in\Xq_0$ and we get a projection $\ol\cU\auf\fa^*$.

\Lemma DDD. For $x\in\Xq$ let $W(x)$ be the little Weyl group of
$\fm_x$. Then $W(x)\subseteq W_X$. Moreover, let $\ol\psi_x$ be the
restriction of $\ol\psi$ to $\Zq_x=\fm_x^\perp$. Then $\ol\psi_x$ is
the composition
$$\eqno{7}
\fm_x^\perp\Pf{\psi_{\fm_x}}\fa^*/W(x)\auf\fa^*/W_X.
$$

\Proof: The compactification $\Zq$ is constructed such that the
restriction of $\ol\mu:\Zq\pfeil\fg^*$ to a fiber $\Zq_x$ is the
natural embedding $\fm_x^\perp\into\fg^*$. Thus we get a commutative
diagram (without dotted arrow)
$$\eqno{}
\matrix{\fm_x^\perp&\into&\Zq\cr
\untenPf&&\untenPf\cr
\fa^*/W(x)&\dashedrightarrow&\fa^*/W_X\cr
&\searrow&\untenPf\cr
&&\fg^*\mod G\cr}
$$
The universal property (as integral closure) implies the existence of
the dotted arrow. In particular, $W(x)\subseteq W_X$.\qed

\Lemma CCC. Let $C$ be an irreducible component of a fiber
$\ol\psi^{-1}(u)$and let $x\in\Xq$. Then $C\cap\Zq_x$ is of pure
codimension $r$ in $\Zq_x$. Moreover, $C$ is of pure codimension $r$
in $\Zq$ and the projection $C\pfeil\Xq$ is dominant.

\Proof: Clearly, the codimension of $C$ in $\Zq$ is at most
$r$. Therefore, the codimension of every irreducible component of
$C\cap\Zq_x$ in $\Zq_x$ is also at most $r$. On the other hand,
$C\cap\Zq_x$ is contained in a fiber of $\ol\psi_x$ which is of pure
codimension $r$ by \cite{equilie}\r2 and \cite{E7}. Hence $C\cap\Zq_x$
is of pure codimension $r$ in $\Zq_x$.

Now consider the projection $p:C\pfeil\Xq$. Let $Y$ be the closure of
its image. We have seen that the fibers of $p$ have dimension
$n-r$. Hence
$$\eqno{}
2n-r\le\|dim|C=(n-r)+\|dim|Y\le(n-r)+n=2n-r
$$
which implies $\|dim|C=2n-r$ and $\|dim|Y=n$, hence $Y=X$.\qed

\Corollary AAA. For $u\in\fa^*/W_X$ consider the fibers
$F=\psi^{-1}(u)\subseteq Z$ and
$\Fq=\ol\psi^{-1}(u)\subseteq\Zq$. Then $\Fq$ is the closure of $F$ in
$\Zq$.

\Proof: By \cite{CCC}, every irreducible component of
$\ol\psi^{-1}(u)$ meets the open subset $\ol\psi^{-1}(u)\cap
Z=\psi^{-1}(u)$. Hence, $\ol\psi^{-1}(u)$ is the closure of $F$.\qed

\beginsection connected. The proof of the connectedness theorem in the
homogeneous case

First, we construct very special points in the fiber $\Fq$ by applying
the following result:

\Proposition. Let $Y$ be an affine $\G_m$\_variety, $V$ a finite
dimensional vector space and $\phi:V\pfeil Y$ a $\G_m$\_equivariant
morphism. Assume that $V$ has a closed irreducible $\G_m$\_stable
subvariety $S$ such that the restriction $\phi|_S:S\pfeil Y$ is finite
and surjective. Let $C$ be an irreducible component of a fiber of
$\phi$. Then $C\cap S\ne\emptyset$.

\Proof: Let $C$ be an irreducible component of $F:=\phi^{-1}(y)$,
$y\in Y$. Assume first $y=y_0:=\phi(0)$. Then $C$ is a closed
$\G_m$\_stable subset of $V$. Hence $0\in C\cap S$ and we are
done. Thus we may assume $y\ne y_0$.

The origin $0$ is the only closed orbit of $V$, hence of $S$. Since
$\phi|_S$ is finite and surjective, the fixed point $y_0$ is the only
closed orbit of $Y$. This implies in particular, that $\G_my$ is
closed in $Y':=Y\setminus\{y_0\}$. Let $V_0:=\phi^{-1}(y_0)$ and
$V':=V\setminus V_0=\phi^{-1}(Y')$. Then $\G_mF=\phi^{-1}(\G_my)$ is
closed in $V'$. Let $E\subset\G_m$ be the
isotropy group of $y$. From $\G_mF\auf\G_my=\G_m/E$ we obtain
$\G_mF=\G_m\times^EF$. The group $E$ is necessarily a finite (cyclic)
group. Hence the map $\G_m\times F\pfeil\G_mF$ is proper, which
implies that $Z:=\G_mC$ is closed in $\G_mF$. We conclude that $Z$ is
also closed in $V'$ or, equivalently,
$$\eqno{1}
\Zq\setminus Z\subseteq V_0.
$$

Let $r:=\|dim|S$. Using again that $\phi|_S$ is finite and surjective
we have that $Y$ is irreducible of dimension $r$. Since $C$ is an
irreducible component of a fiber of $\phi:V\pfeil Y$, we have
$\|codim|_VC\le\|dim|Y=r$. Thus, $\|codim|_V\Zq<r$. Since $0\in\Zq$
and $\|dim|S=r$, the intersection $\Zq\cap S$ is non\_empty of
positive dimension (here we use the smoothness of $V$). On the other
hand, the finiteness of $\phi|_S$ implies that $V_0\cap S$ is a finite
set. Thus we get from \cite{E1} that $Z\cap S\ne\emptyset$ which is
equivalent to $C\cap S\ne\emptyset$.\qed

\Corollary intersection. Let $\fa'\subseteq\fh^\perp$ as in
\cite{equilie}\r4 and let $C$ be an irreducible component of a
fiber of $\psi_\fh:\fh^\perp\pfeil\fa^*/W_X$. Then
$C\cap\fa'\ne\emptyset$.

\noindent Next we study the intersection $C\cap\fa'$ for the Lie
algebras $\fm_x$ simultaneously for an open set of $x\in X_0$. For
this, let $\cA\subseteq\ol\cU$ be a complementary vector bundle to
$\ol\cB$ over an open subset $\Xq_1\subseteq\Xq_0$. This exists in the
neighborhood of every point of $\Xq_0$. Since the fiber is
$\cA_x=\ol\cU_x/\ol\cB_x=\fa^*$, there is a canonical trivialization
$\tau:\fa^*\times\Xq_1\pf\sim\cA$.

\Lemma BBB. Let $C$ be an irreducible component of the fiber
$F:=\ol\psi^{-1}(u)$ and assume there is $\alpha\in\fa^*$, $x_0\in \Xq_1$
with $\tau(\alpha,x_0)\in C$. Then $\tau(\alpha,x)\in C$ for every
$x\in\Xq_1$.

\Proof: The restriction of $\ol\psi$ to $\cA$ is the composition
$\cA\auf\fa^*\pf\pi\fa^*/W_X$. Therefore, $F\cap\cA$ is of pure
dimension $n=\|dim|X$. More precisely, $F\cap \cA=\tau(S\times\Xq_1)$
where $S$ is the $W_X$\_orbit $\pi^{-1}(u)$. On the other hand, $C$
is of pure codimension $r$ in $\Zq$ and $\cA$ is irreducible of
dimension $n+r$. This implies that the dimension of every irreducible
component of $C\cap\cA$ is at least $n$. Since $C\cap\cA\subseteq
F\cap\cA$, we conclude that $C\cap\cA$ is the union of irreducible
components of $F\cap\cA\cong S\times\Xq_1$. Hence there is a subset
$S'$ of $S$ with $C\cap\cA=\tau(S'\times\Xq_1)$. By assumption we have
$(\alpha,x_0)\in S'\times\Xq_1$. Hence $\alpha\in S'$ and
$\tau(\alpha,x)\in C\cap\cA\subseteq C$ for all $x\in\Xq_1$.\qed

Now we are able to prove \cite{ConnThm} for $X=G/H$. Consider the
fiber $F=\psi^{-1}(u)$ and suppose that $F$ is disconnected. Then
$F=F_1\dot\cup F_2$ is the disjoint union of non\_empty closed
subsets.  By \cite{AAA}, the closure of $F$ in $\Zq$ is the fiber
$\Fq=\ol\psi^{-1}(u)$ which is connected by \cite{EEE}. This implies
$\Fq_1\cap\Fq_2\ne\emptyset$ but we can say even more. As fiber of an
equidimensional map between smooth varieties, $F$ is locally a
complete intersection. A theorem of Hartshorne (\cite{Ha}~3.4) then
asserts that $F$ is even connected in codimension one, i.e., stays
connected upon removal of any subset of codimension two or
higher. Applied to our situation, we conclude that $\Fq_1\cap\Fq_2$
has an irreducible component $C_0$ which is of codimension $1$ in
$\Fq$, hence of codimension $r+1$ in $\Zq$.

Since $C_0\cap Z\subseteq \Fq_1\cap\Fq_2\cap Z=F_1\cap F_2=\emptyset$
we have $C_0\cap Z=\emptyset$. Let $X'$ be the closure of the image of
$C_0$ in $\Xq$. Then $X'\cap X=\emptyset$. Let $x\in X'$. Then
$C_0\cap\Zq_x$ is a subset of a fiber of $\ol\psi_x$, hence has
dimension at most $n-r$ (\cite{CCC}). This implies
$$\eqno{}
2n-r-1=\|dim|C_0\le (n-r)+\|dim|X'\le (n-r)+(n-1)=2n-r-1. 
$$ Since equality has to hold throughout, we get in particular that
there exists an $x_0\in X'$ such that
$\|dim|C_0\cap\Zq_{x_0}=n-r$. Therefore, $C_0\cap\Zq_{x_0}$ is the
union of irreducible components of a fiber of $\ol\psi_x$. Now choose
a complement $\cA$ of $\ol\cB\subseteq\ol\cU$ in a neighborhood of
$x_0$. From \cite{intersection}, we get
$C_0\cap\cA_{x_0}\ne\emptyset$. In particular, there is
$\alpha\in\fa^*$ with $\tau(\alpha,x_0)\in C_0$.

Now let $C_1\subseteq F_1$, $C_2\subseteq F_2$ be irreducible
components such that $C_0\subseteq\Cq_1\cap\Cq_2$. Since
$\tau(\alpha,x_0)\in\Cq_i$ we have $\tau(\alpha,x)\in\Cq_i$ for all
$x$ in a neighborhood of $x_0$ (\cite{BBB}). If we choose $x\in X$ we
get $\tau(\alpha,x)\in C_1\cap C_2\subseteq F_1\cap F_2$ in
contradiction to $F_1\cap F_2=\emptyset$. This finishes the proof of
\cite{ConnThm} when $X$ is homogeneous.

\beginsection generalcase. The proof in general

From now on, let $X$ be any smooth $G$\_variety.  To prove the general
case we use the following trivial lemma:

\Lemma triviallemma. Let $F$ be a variety and $F'\subseteq F$ an open
subset. Assume
\item{\r1}$F'$ is connected and
\item{\r2}for every point $\eta\in F$ there are two irreducible
subvarieties $D$ and $E$ of $F$ such that
$\eta\in D$, $D\cap E\ne\emptyset$, and $E\cap F'\ne\emptyset$.
\Par\noindent
Then $F$ is connected.\qed

\noindent Fix $u\in\fa^*/W_X$ and let $F:=\psi^{-1}(u)$. To construct
$F'$ we start with:

\Lemma. Let $X$ be an irreducible $G$\_variety. Then there is a
non\_empty open $G$\_stable subset $X'\subseteq X$ with the following
properties:
\item{\r1}The orbit space $Q=X'/G$ exists, i.e., there is a
$G$\_invariant surjective morphism $\pi:X'\pfeil Q$ such that all
fibers $X'_v:=\pi^{-1}(v)$ are reduced and homogeneous.
\item{\r2}The spaces $X'$ and $Q$ are smooth.
\item{\r3}All orbits $X'_v$ share the same dimension, complexity, rank,
and little Weyl group $W_X$.\Par

\Proof: The existence of an open subset $X'$ as in \r1 is a well known
theorem of Rosenlicht (\cite{Ros}, see also \cite{Sp}~\S2). Then \r2
can be achieved by shrinking $Q$ to an open subset. Finally, \r3
follows from \cite{WuM}~2.5 (for dimension, complexity, and rank) and
\cite{WuM}~6.5.4 (for the little Weyl group).\qed

\noindent Let $X'\subseteq X$ be as in the lemma. Let
$$\eqno{}
T^*_{X'/Q}:=\dotunion_{v\in Q}T^*_{X'_v}\pfeil X'
$$
be the relative cotangent bundle. Since $W_{X'_v}=W_X$ we get
morphisms $\psi_v:T^*_{X'_v}\pfeil\fa^*/W_X$ which glue to a morphism
$\psi_*:T^*_{X'/Q}\pfeil\fa^*/W_X$. Let $N\subseteq T^*_{X'/Q}$ be the
fiber $\psi^{-1}_*(u)$ and let $\pi:N\pfeil Q$ be the projection. Then
$\pi^{-1}(v)=\psi^{-1}_v(u)$ which is connected since $X'_v$ is
homogeneous. Upon shrinking $Q$ to an open subset we may assume that
$\pi$ is flat and surjective. Thus, \cite{openclosed} implies that $N$
is connected. Next observe that there is a projection $p:T^*_{X'}\auf
T^*_{X'/Q}$. In fact, this is a locally trivial bundle of affine
spaces and therefore flat, surjective with connected fibers. On the
other hand, the morphism $\psi:T^*_{X'}\pfeil\fa^*/W_X$ is just the
composition $\psi_*\circ p$. Thus
$$
F':=F\cap T^*_{X'}=p^{-1}(N)
$$
is connected which verifies part \r1 of \cite{triviallemma}.

For the second part, let $\eta$ be any point of $F$ and let $Y$ be the image
of $G\eta$ in $X$. We need to compare the moment map for $X$ with that
for the orbit $Y\subseteq X$. The main tool to study the
interrelation is the local structure theorem:

\Theorem. There exists a parabolic subgroup $P$ (containing $B$) with
Levi part $L$ (containing $T$) and an affine $L$\_stable subvariety
$S$ of $X$ with:
\item{\r1}The natural morphism $P\times^LS\pfeil X$ is an open embedding.
\item{\r2}The intersection $Y_0:=Y\cap S$ is a non\_empty $L$\_variety.
\item{\r3}Let $L_0$ be the kernel of the action of $L$ on $Y_0$. Then
$A_Y:=L/L_0$ is a torus.
\item{\r4}The action of $A_Y$ on $Y_0$ is locally free, i.e., has
finite isotropy groups.\Par

\def\loccit{{\it loccit.}}

\Proof: This theorem is essentially due to Brion\_Luna\_Vust
\cite{BLV}. We use the refinement in \S\S1--2 of \cite{IB} as
follows. By \loccit~2.1, we may replace $X$ by a $G$\_stable open
subset which supports an ample $G$\_linearized line bundle $\cL$. Now
we choose a section $\sigma$ of a power of $\cL$ as in
\loccit~2.10. Let $P=G_{\<\sigma\>}$ be the stabilizer of the line
$\<\sigma\>=\CC\sigma$.

Now \r1 follows from \loccit~1.2.3 while \r2 is implied by
\loccit~2.10.2.  Furthermore, \loccit~2.10.3 and 2.8.1 together show
\r3. Finally, by \loccit~2.10.4 all orbits of $A_Y$ in $Y_0$ are
closed. Since $Y_0$ is affine this implies that they have all the same
dimension. Since the generic isotropy group is trivial by construction
we see that all isotropy groups are finite.\qed

\noindent Let $P_u$ be the unipotent radical of $P$. Since
$P=P_u\times L$ (as a right $L$\_variety) we have $P\times^LS=P_u\times
S$. Thus,
$$\eqno{}
P_u\times S\pfeil X
$$
is an open embedding. This implies that $S$ is smooth and
irreducible. Moreover, for every $x\in S$ the tangent space splits
$$\eqno{}
T_{X,x}=\fp_ux\oplus T_{S,x}.
$$
This allows us to embed $T^*_S$ into $T^*_X$. The relationship between
the moment maps is given by:

\Lemma. The image of $T^*_S$ under the moment map on $T^*_X$ is in
$\fp_u^\perp$. Moreover, the following diagram commutes:
$$\eqno{9}
\matrix{
T^*_S&&\hbox to 0pt{\hss$\into$\hss}&&T^*_X\cr
\Links{\mu_S}\untenPf&\searrow&&&\Links{\mu_X}\untenPf\cr
\fl^*&\leftarrow&\!\!\!\fp_u^\perp\!\!\!&\into&\fg^*\cr
\untenPf&&&&\untenPf\cr
\fl^*\mod L&&\hbox to 0pt{\hss$\Pfeil$\hss}&&\fg^*\mod G\cr}
$$

\Proof: That image of $T^*_S$ under the moment map on $T^*_X$ is, in
$\fp_u^\perp$, is a reformulation of the definition of the embedding
$T^*_S\into T^*_X$. This proves the commutativity of the
quadrangle. The triangle commutes since $S\into X$ is
$L$\_equivariant. Finally, the commutativity of the pentagon is
\cite{WuM} 6.1\qed

Now choose $y\in Y_0$. Since $S$ is affine and since the orbit $Ly$ is
closed, it has a ``slice'' in $S$, i.e., an irreducible, smooth,
$L_0$\_stable subvariety $S_0$ of $S$ with $y\in S_0$ and
such that $L\times^{L_0}S_0\pfeil S$ is \'etale. The Lie algebra of
$L$ decomposes as $\fl=\fa_Y\oplus\fl_0$ which induces a decomposition of tangent spaces: $T_{S,x}=\fa_Yx\oplus
T_{S_0,x}$ for every $x\in
S_0$. Thus, we get embeddings and commutative diagrams
$$\eqno{10}
\matrix{\fa_Y^*\times S_0&\into& T^*_S&\into& T^*_X\cr
\untenPf&&\untenPf&&\untenPf\cr
\fa_Y^*&\pfeil&\fl^*\mod L&\pfeil&\fg^*\mod G.\cr}
$$
Now we look more closely at the subspace
$\fa^0:=\fa_Y^*\times\{y\}\subseteq T^*_{X,y}$. We have
$$\eqno{8}
T_{X,y}=\fp_uy\oplus\fa_Yy\oplus T_{S_0,y}
$$
Since $\fa_Yy=\fl y$ we have
$$\eqno{}
\fp_uy\oplus\fa_Yy=\fb y\quad\hbox{and}\quad\fp_uy=\fu y.
$$
Thus, $\fa^0$ is a complement of $(\fb y)^\perp$ in $(\fu y)^\perp$.

On the other hand, since $P\times^LY_0\pfeil Y$ is an open embedding,
the generic $P$\_orbit in $Y$ is isomorphic
to $P/L_0$. Because $Py=P/L_0=B/(B\cap L_0)=By$ it follows that $By$
represents a generic $B$\_orbit in $Y$. Moreover, $Uy$ is a generic
$U$\_orbit in $Y$.

Now consider the restriction map $\rho:T^*_{X,y}\auf
T^*_{Y,y}=\fg_y^\perp$. The moment map $\mu$ clearly factors through
$\rho$. Moreover, the image $\fa':=\rho(\fa^0)$ is a complement of
$(\fb y)^\perp$ in $(\fu y)^\perp$. It follows from
\cite{intersection} that $\fa'$ intersects every irreducible component
of every fiber of $T^*_{Y,y}\pfeil\fg^*\mod G$. Thus the same holds
for $\fa^0$ with respect to the morphism $T^*_{X,y}\pfeil\fg^*\mod G$.

Now we are able to prove the connectedness theorem in general: by
changing $\eta$ in its $G$\_orbit we may assume that $\eta\in
T^*_{X,y}\cap F$. Let $D$ be an irreducible component of the fiber of
$T^*_{X,y}\pfeil\fg^*\mod G$ containing $\eta$. By construction, $F$
is a union of connected components of a fiber of $T^*_X\pfeil\fg^*\mod
G$. Thus $D\subseteq F$. Moreover, there is $\alpha\in\fa^*_Y$ with
$(\alpha,y)\in D$. Let $E:=\{\alpha\}\times S_0$. Then the commutative
diagrams \cite{E9} and \cite{E10} show that $E$ is contained in a
fiber of $T^*_X\pfeil\fg^*\mod G$. Thus $E\subseteq F$ and $D\cap
E\ne\emptyset$. Finally, $G\cdot S_0$ is dense in $X$. Hence $S_0\cap
X'\ne\emptyset$ which implies $E\cap F'\ne\emptyset$. Thus, we
conclude with \cite{triviallemma} that $F$ is connected.

\beginrefs

\L|Abk:An|Sig:An|Au:Anantharaman, S.|Tit:Sch\'emas en groupes, espaces
homog\`enes et espaces alg\'ebriques sur une base de
dimension 1|Zs:Bull. Soc. Math. France, Mem.|Bd:33|S:5--79|J:1973|xxx:-||

\B|Abk:AHV|Sig:AHV|Au:Aroca, J.; Hironaka, H.; Vicente, J.%
|Tit:Desingularization theorems|Reihe:Mem\'orias de Matem\'atica del Instituto
``Jorge Juan'' {\bf 30}|Verlag:Consejo Superior de Investigaciones
Cient\'\i ficas|Ort:Madrid|J:1977|xxx:-||

\L|Abk:BLV|Sig:BLV|Au:Brion, M.; Luna, D.; Vust, Th.|Tit:Espaces homog\`enes
sph\'eriques|Zs:Invent. Math.|Bd:84|S:617--632|J:1986|xxx:-||

\B|Abk:DG|Sig:DG|Au:Demazure, M.; Gabriel, P.|Tit:Groupes alg\'ebriques.
Tome I|Reihe:-|Verlag:North-Holland Publishing Co.|Ort:Amsterdam|J:1970|xxx:-||

\L|Abk:EGA|Sig:EGA|Au:Dieudonn\'e, J.; Grothendieck, A.|Tit:El\'ements de
g\'eom\'etrie alg\'ebrique IV|Zs:Publ. Math. IHES|Bd:28|S:-|J:1966|xxx:-||


\L|Abk:Ha|Sig:\\Ha|Au:Hartshorne, R.|Tit:Complete intersections and
connectedness|Zs:Amer. J. Math.|Bd:84|S:497--508|J:1962|xxx:-||

\B|Abk:Har|Sig:\\Ha|Au:Hartshorne, R.|Tit:Algebraic
Geometry|Reihe:Graduate Texts in
Mathematics~{\bf 52}|Verlag:Springer|Ort:Heidelberg|J:1977|xxx:-||

\L|Abk:Kir|Sig:Ki|Au:Kirwan, F.|Tit:Convexity properties of the moment
mapping III|Zs:Invent. math.|Bd:77|S:547-552|J:1984|xxx:-||

\L|Abk:WuM|Sig:\\Kn|Au:Knop, F.|Tit:Weylgruppe und Momentabbildung|%
Zs:Invent. Math.|Bd:99|S:1-23|J:1990|xxx:-||

\L|Abk:IB|Sig:\\Kn|Au:Knop, F.|Tit:{\"U}ber Bewertungen, welche unter einer
reduktiven Gruppe invariant sind|Zs:Math. Ann.|Bd:295|S:333--363|J:1993|xxx:-||

\L|Abk:Motion|Sig:\\Kn|Au:Knop, F.|Tit:The asymptotic behavior of invariant
collective motion|Zs:Invent. math.|Bd:116|S:309--328|J:1994|xxx:-||



\L|Abk:WHM|Sig:\\Kn|Au:Knop, F.|Tit:Weyl groups of Hamiltonian manifolds,
I|Zs:Preprint|Bd:-|S:34 pages|J:1997|xxx:dg-ga/9712010||

\L|Abk:Ros|Sig:Ro|Au:Rosenlicht, M.|Tit:On quotient varieties and the affine
embedding of certain homogeneous spaces|Zs:Trans. Am. Math.
Soc.|Bd:101|S:211--223|J:1961|xxx:-||

\Pr|Abk:Sp|Sig:Sp|Au:Springer, T.|Artikel:Aktionen reduktiver Gruppen auf
Variet\"aten|Titel:Algebraische Transformationsgruppen und Invariantentheorie|Hgr:H. Kraft, P. Slodowy, T. Springer, eds.|Reihe:DMV Sem.|Bd:13|Verlag:Birkh{\"a}user|Ort:Basel|S:3--39|J:1989|xxx:-||


\endrefs

\bye